\documentclass[10pt]{article}
\usepackage{amsthm}
\usepackage{latexsym}
\usepackage{amsmath}
\usepackage{amssymb}
\usepackage{epsfig}
\usepackage[matrix,arrow,cmtip,curve,dvips]{xy}

    \newtheorem{thm}{Theorem}[section]
    \newtheorem{prop}[thm]{Proposition}
    \newtheorem{lemma}[thm]{Lemma}
    
    \newtheorem{cor}[thm]{Corollary}

    \newtheorem{defn}[thm]{Definition}

    \newtheorem{rem}[thm]{Remark}
    \newtheorem{example}[thm]{Example}

\newcommand{\Hom}{\mbox{\upshape Hom}}

\newcommand{\colim}{\mbox{\upshape colim}}

\newcommand{\im}{\mbox{\upshape im}}

\title{Homotopy groups of $\Hom$ complexes of graphs}

\author{Anton Dochtermann \\
\small Institut f\"{u}r Mathematik, MA 6-2 \\[-0.8ex]
\small Technische Universit\"{a}t Berlin \\[-0.8ex]
\small Strasse des 17. Juni 136 \\[-0.8ex]
\small 10623 Berlin, Germany \\[-0.8ex]
\small \texttt{anton.dochtermann@gmail.com} }

\begin{document}

\maketitle

\begin{abstract}
The notion of $\times$-homotopy from \cite{DocHom} is investigated
in the context of the category of pointed graphs.  The main result
is a long exact sequence that relates the higher homotopy groups of
the space $\Hom_*(G,H)$ with the homotopy groups of $\Hom_*(G,H^I)$.  Here $\Hom_*(G,H)$ is a space which parameterizes pointed graph maps from $G$ to $H$ (a pointed version of the usual $\Hom$ complex), and $H^I$ is the graph of based paths in $H$.
As a corollary it is shown that $\pi_i \big(\Hom_*(G,H) \big) \cong [G,
\Omega^i H]_{\times}$, where $\Omega H$ is the graph of based closed paths
in $H$ and $[G,K]_{\times}$ is the set of $\times$-homotopy classes of pointed graph maps from $G$ to $K$. This is similar in spirit to the results of \cite{BBLL},
where the authors seek a space whose homotopy groups encode a
similarly defined homotopy theory for graphs.  The categorical
connections to those constructions are discussed.
\end{abstract}

\section{Introduction}
In several recent papers (see for instance \cite{BBLL}, \cite{BL}), a homotopy
theory of reflexive graphs termed $A$-theory has been developed in
which graph theoretic homotopy groups are defined to measure
`combinatorial holes' in simplicial complexes. In \cite{BBLL} the
authors construct a cubical complex $X_G$ (associated to the graph
$G$), and a homomorphism from the homotopy groups of the geometric
realization of $X_G$ to the $A$-theory groups of $G$; modulo a (yet
unproved) version of cubical approximation they show that this map
is in fact an isomorphism.

In the paper \cite{DocHom}, a similar homotopy theory for general
graphs called $\times$-homotopy is developed.  Both theories are
discussed in the common framework of exponential graph constructions
associated to the relevant product (cartesian for $A$-theory,
categorical for $\times$-homotopy).  There it is shown that
$\times$-homotopy is characterized by topological properties of the
so-called $\Hom$ complex, a functorial way to assign a poset (and
hence topological space) to a pair of graphs, first introduced to
provide lower bounds on the chromatic number of graphs.  In
particular, the $\times$-homotopy class of maps from graphs $G$ to
$H$ are seen to coincide with the path components of the space
$\Hom(G,H)$.

In this paper, we consider the graph theoretic notions of homotopy
groups that arise in the context of $\times$-homotopy.  In order to give a topological interpretation of these constructions it is necessary to restrict our attention to the category of pointed graphs.  We show that
these combinatorially defined groups are isomorphic to the usual homotopy groups of a
pointed version of the $\Hom$ complex (denoted $\Hom_*$).  Our
method for proving this is to construct a `path graph' $G^I$
associated to a (pointed) graph $G$ and to show that the natural
endpoint map induces a long exact sequence of homotopy groups of the
$\Hom_*$ complexes.

The paper is organized as follows.  In Section 2, we describe the
category of pointed graphs and recall the notions of both $A$-homotopy and
$\times$-homotopy.  In Section 3, we introduce the $\Hom_*$ functors
and establish some basic facts about them, including interaction with
the relevant adjunction and also a graph operation known as folding.
In Section 4, we introduce the notion of a path graph, and state and
prove our main result regarding the long exact sequence of the
homotopy groups of the $\Hom_*$ complexes.  In Section 5, we end with a brief discussion regarding other aspects of our construction.

\noindent \textbf{Acknowledgments.}  The author wishes to thank his advisor Eric Babson for guiding this research, and also the two anonymous referees for providing valuable suggestions which helped to improve the exposition.

\section{Basic objects of study}
We begin with the basic definitions. For us, a \textit{graph} $G =
\big(V(G), E(G)\big)$ consists of a vertex set $V(G)$ and an edge set $E(G)
\subseteq V(G) \times V(G)$ such that if $(v,w) \in E(G)$ then
$(w,v) \in E(G)$. Hence our graphs are undirected and do not have
multiple edges, but may have loops (if $(v,v) \in E(G)$).  If $(v,w)
\in E(G)$ we will say that $v$ and $w$ are \textit{adjacent} and
denote this as $v \sim w$. Given a pair of graphs $G$ and $H$, a
\textit{graph homomorphism} (or \textit{graph map}) is a mapping of the
vertex set $f:V(G) \rightarrow V(H)$ that preserves adjacency: if $v
\sim w$ in $G$, then $f(v) \sim f(w)$ in $H$.  With these as our
objects and morphisms we obtain a category of graphs which we denote
${\mathcal G}$.  If $v$ is a vertex of a graph $G$, then $N(v) :=
\{w \in V(G): (v,w) \in E(G)\}$ is called the \textit{neighborhood} of $v$.  A \textit{reflexive} graph is a graph $G$ with
loops on all the vertices, so that $v \sim v$ for all $v \in V(G)$.
The category of reflexive graphs is denoted ${\mathcal
G}^\circ$.  For more about graph homomorphisms and the category of
graphs, see \cite{HN04}.

There are two relevant monoidal structures on the category of
graphs.  The first is the \textit{categorical product} $G \times H$
of graphs $G$ and $H$, defined to be the graph with vertex set $V(G)
\times V(H)$ and adjacency given by $(g,h) \sim (g^\prime,h^\prime)$
if $g \sim g^\prime$ and $h \sim h^\prime$.  The other is the
\textit{cartesian product} $G \square H$, defined to be the graph
with the same vertex set $V(G) \times V(H)$ with adjacency given by
$(g,h) \sim (g^\prime,h^\prime)$ if either $g \sim g^\prime$ and $h
= h^\prime$, or $h \sim h^\prime$ and $g = g^\prime$.  Each of these
products has a right adjoint given by versions of an internal hom
construction (see \cite{DocHom} for further discussion).

In this paper, we work primarily in the category of pointed graphs.
A \textit{pointed graph} $G = (G,x)$ is a graph $G$ together with a
specified looped vertex $x$.  A \textit{map of pointed graphs}
$f:(G,x) \rightarrow (H,y)$ is a map of graphs such that $f(x) = y$.
The resulting category of pointed graphs will be denoted ${\mathcal
G}_*$.  The category ${\mathcal G}_*$ enjoys some useful properties
which we discuss next.

\begin{defn}
For pointed graphs $G = (G,x)$ and $H = (H,y)$, the \textit{smash
product} $G \wedge H$ is the pointed graph with vertex set given by
the quotient of $(V(G) \times V(H))$ under the identifications
$(x,h) = (g,y)$ for all $g \in V(G)$, $h \in V(H)$. Adjacency is
given by $[(g,h)] \sim [(g^\prime, h^\prime)]$ if $g \sim g^\prime$
and $h \sim h^\prime$ for some representatives. The graph $G \wedge
H$ is pointed by the vertex $[(x,y)]$ (see Figure 1).
\end{defn}

\vspace{.1 in}
\begin{center}

\epsfig{file=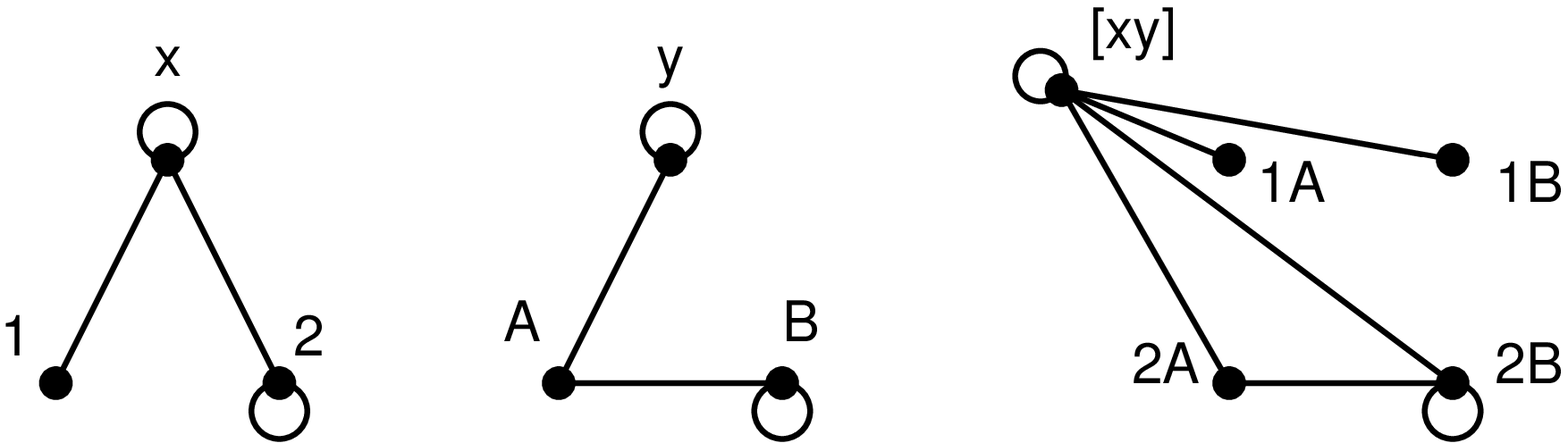, height=1 in, width = 3.5 in}

{Figure 1: The graphs $G=(G,x)$, $H=(H,y)$, and $G \wedge H = (G
\wedge H,[(x,y)])$}
\end{center}

\begin{defn}For pointed graphs $G = (G,x)$ and $H = (H,y)$, the
\textit{pointed} internal hom graph, denoted $H^G$, is the pointed
graph with vertices given by all \textit{set} maps $\{f:V(G)
\rightarrow V(H): f(x) = y\}$.  Adjacency is given by $f \sim g$ if
$f(v) \sim g(v^\prime)$ in $H$ for all $v \sim v^\prime$ in $G$. The
graph $H^G$ is pointed by the graph map that sends every vertex of
$G$ to the vertex $y \in H$.
\end{defn}

These last two constructions are adjoint to one another, as
described by the following lemma.  We let ${\mathcal G}_*(G,H)$
denote the set of pointed maps between the pointed graphs $G$ and
$H$.

\begin{lemma}\label{pointadjunct}
For pointed graphs $A = (A,x)$, $B = (B,y)$, and $C = (C,z)$ we have
a natural bijection of sets $\varphi:{\mathcal G}_*(A \wedge B, C)
\rightarrow {\mathcal G}_*(A, C^B)$, given by $\varphi(f)(a)(b) =
f[(a,b)]$ for $a \in V(A)$ and $b \in V(B)$.
\end{lemma}

\begin{proof}
First note that $\varphi(f)$ is well-defined since $f(x,b) = f(a,y)
= z$ for all $a \in V(A)$, $b \in V(B)$.  Next we see that
$\varphi(f)(a) \in C^B$ since $\varphi(f)(a)(y) = f(a,y) = z$ and $\varphi(f)(x)(b) = f(x,b) = z$.  To see that $\varphi(f)$ is a graph map, suppose $a \sim
a^\prime$; we need to check that $\varphi(f)(a) \sim
\varphi(f)(a^\prime)$.  If $b \sim b^\prime$ then we have
$\varphi(f)(a)(b) = f(a,b)$ and $\varphi(f)(a^\prime)(b^\prime) =
f(a^\prime, b^\prime)$, and hence (the equivalence classes) are
adjacent in $A \wedge B$.

Next, to show that $\varphi$ is a bijection we define a map
$\psi:{\mathcal G}_*(A,C^B) \rightarrow {\mathcal G}_*(A \wedge B,
C)$ via $\psi(g)[(a,b)] = g(a)(b)$.  We note that $\psi(g)(x,b) =
g(x)(b) = z$ and $\psi(g)(a,y) = g(a)(y) = z$ for $a \in V(A)$ and
$b \in V(B)$, and hence $\psi(g)$ is well defined on the vertices of
$A \wedge B$.  Similarly, one checks that $\psi(g)$ is a pointed
graph map.  It is clear that $\psi$ is the inverse to $\varphi$ and
the result follows.
\end{proof}

We next turn to the definition of `homotopy' in the context of the pointed category.  The construction runs along the same lines of
the $\times$-homotopy introduced in \cite{DocHom}.  Although the
constructions for the pointed category are straightforward
modifications of the notions for the unpointed category, we will record all the definitions here for convenience.

If $G$ is a graph we define $G_*$ to be the pointed graph
obtained by adding a distinguished disjoint looped vertex (denoted $*$) to the
graph $G$.

\begin{defn}
The graph $I_n$ is the (reflexive) graph with vertices $\{0,1,
\dots, n\}$ and with adjacency given by $i \sim j$ if $|i-j| \leq
1$ (see Figure 2).
\end{defn}

\begin{center}
\epsfig{file=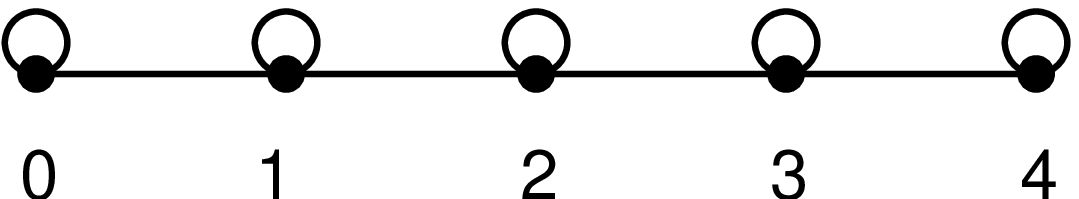, height=.3 in, width = 1.8 in}

{Figure 2: The graph $I_4$}
\end{center}

For our definition of homotopy, we will want to consider path
components of the exponential graph $H^G$.  To capture this notion
in the category of pointed graphs, the relevant graph to map from
will be the pointed graph $I_{n *}$, since a pointed map $I_{n *}
\rightarrow (G,x)$ will send $*$ to $x$, and $0$ and $n$ to the
endpoints of a path in $G$.

\begin{defn}
A pair of maps $f,g: G \rightarrow H$ between pointed graphs is
called \textit{$\times$-homotopic} if there is an integer $n$ and a
(pointed) graph map $F:I_{n *} \rightarrow H^G$ such that $F(0) = f$
and $F(n) = g$.
\end{defn}

This defines an equivalence relation on ${\mathcal G}_*(G,H)$, and
the set of $\times$-homotopy classes of pointed maps between $G$ and
$H$ will be denoted $[G,H]_{\times}$ (or simply $[G,H]$ if the context is clear).  Applying the adjunction of
Lemma \ref{pointadjunct} we note that a homotopy between $f$ and $g$
is the same as a pointed map $\tilde F:G \wedge {I_n}_* \rightarrow
H$ with $\tilde F(?,0) = f$ and $\tilde F(?,n) = g$.

In \cite{DocHom}, we relate the construction of our
$\times$-homotopy to the $A$-theory of \cite{BBLL}.  We briefly
recall the discussion here.

\begin{defn}
A pair of pointed graph maps $f,g: (G,x) \rightarrow (H,y)$ between
reflexive graphs is called \textit{$A$-homotopic}, denoted $f \simeq_A g$, if there is an integer $n$ and a graph map

\begin{center}
$\phi: G \square I_n \rightarrow H$,
\end{center}
such that $\phi(?,0) = f$ and $\phi(?,n) = g$, and such that
$\phi(x,i) = y$ for all $i$.

\end{defn}

One can check that $f$ and $g$ are $A$-homotopic if and only if
there is a pointed graph map $\tilde \phi:I_{n *} \rightarrow H^G$
such that $\tilde \phi(0) = f$ and $\tilde \phi (n) = g$; here $H^G$
is the internal hom graph that is right adjoint to the \textit{cartesian}
product.  In other words, our definition of $\times$-homotopy is the
`same' as that of $A$-homotopy, with the categorical product playing
the role of the cartesian product.

In the paper \cite{BBLL}, the authors consider the so-called
\textit{A-homotopy groups} of a reflexive pointed graph $(G,x)$, which they denote $A_n(G,x)$.  By definition, $A_n(G,x)$ is the set of
$A$-homotopy classes of graph maps

\begin{center}
$f:(I^m, \partial I^m) \rightarrow (G,x)$.
\end{center}

\noindent
Here $I^m = I_n \square I_n \square \dots \square I_n$ is the
$m$-fold cartesian product of $I_n$, and $\partial I^m$ is the
subgraph of $I^m$ consisting of vertices with at least one
coordinate equal to $0$ or $n$.

The authors of \cite{BBLL} seek to construct a topological space
whose homotopy groups encode the $A$-homotopy groups of the graph
$G$; this is seen as a generalization of a result from \cite{BL},
where it is shown that $A_1(G,x)$ is isomorphic to the fundamental
group of the space obtained by attaching 2-cells to all 3- and
4-cycles of the graph $G$.  The main result of \cite{BBLL} is the
construction of a cubical complex $M_*(G)$ associated to the
reflexive graph $G$, and a homomorphism from the geometric
realization $X_G := |M_*(G)|$ to the $A$-homotopy groups of $G$.
Here the $i$-cube $M_i(G)$ is defined to be the set ${\mathcal G}(I_1^m, G)$ of
all graph maps from $I_1^m$ to $G$.

In this paper we consider the analogous questions in the context of
$\times$-homotopy.  One can follow the procedure of \cite{BBLL} and
construct a cubical complex built from the sets ${\mathcal G}(I_1^m,
G)$, where this time $I_1^m$ denotes the $m$-fold
\textit{categorical} product.  In this way one obtains a map from the
realization of this space to the graph-theoretically defined
`$\times$-homotopy groups' of $G$.  However, it turns out that we
can follow a somewhat different route to obtain a space
$\Hom_*(T,G)$ whose homotopy groups do in fact coincide with what we
will call the `$T$-homotopy groups of $G$'.  We turn to a discussion
of these spaces in the next section.

\section{$\Hom_*$ complexes}

The $\Hom$ complex is a functorial way to assign a poset (and hence
topological space) to a pair of graphs.  Spaces of this sort were
first introduced by Lov\'{a}sz in \cite{Lov78}, and have more
recently been studied by various people in a variety of contexts (see for example \cite{BKcom}, \cite{BKpro}, \cite{DocUni}, \cite{Ka},
\cite{Kchr}, \cite{MZ}, \cite{Schspace}, and \cite{Ziv05}). The connection to $\times$-homotopy
of ordinary (unpointed) graphs is explored in the paper
\cite{DocHom}. There is a natural notion of the $\Hom$ complex in
the pointed setting (which we will denote as $\Hom_*$). As in the
unpointed setting, this construction interacts well with the
adjunction of Lemma \ref{pointadjunct}, and the path components of
this space characterize the set of $\times$-homotopy of pointed
maps.  We collect these facts next.

\begin{defn}
For pointed graphs $G=(G,x), H=(H,y)$, we define $\Hom_*(G,H)
\subseteq \Hom(G,H)$ to be the (pointed) poset whose elements are
given by all functions $\eta: V(G) \rightarrow 2^{V(H)} \backslash
\{\emptyset\}$, such that $\eta(x) = \{y\}$, and if $(v,w) \in
E(G)$ then for all $\tilde v \in \eta(v)$ and $\tilde w \in
\eta(w)$ we have $(\tilde v, \tilde w) \in E(H)$.  The relation is
given by containment, so that $\eta \leq \eta^\prime$ if $\eta(v)
\subseteq \eta^\prime(v)$ for all $v \in V(G)$.
\end{defn}

\begin{example}
As an example, we consider the pointed graphs $G=(G,x)$ and $H=(H,y)$ in Figure 3.

\begin{center}
\epsfig{file=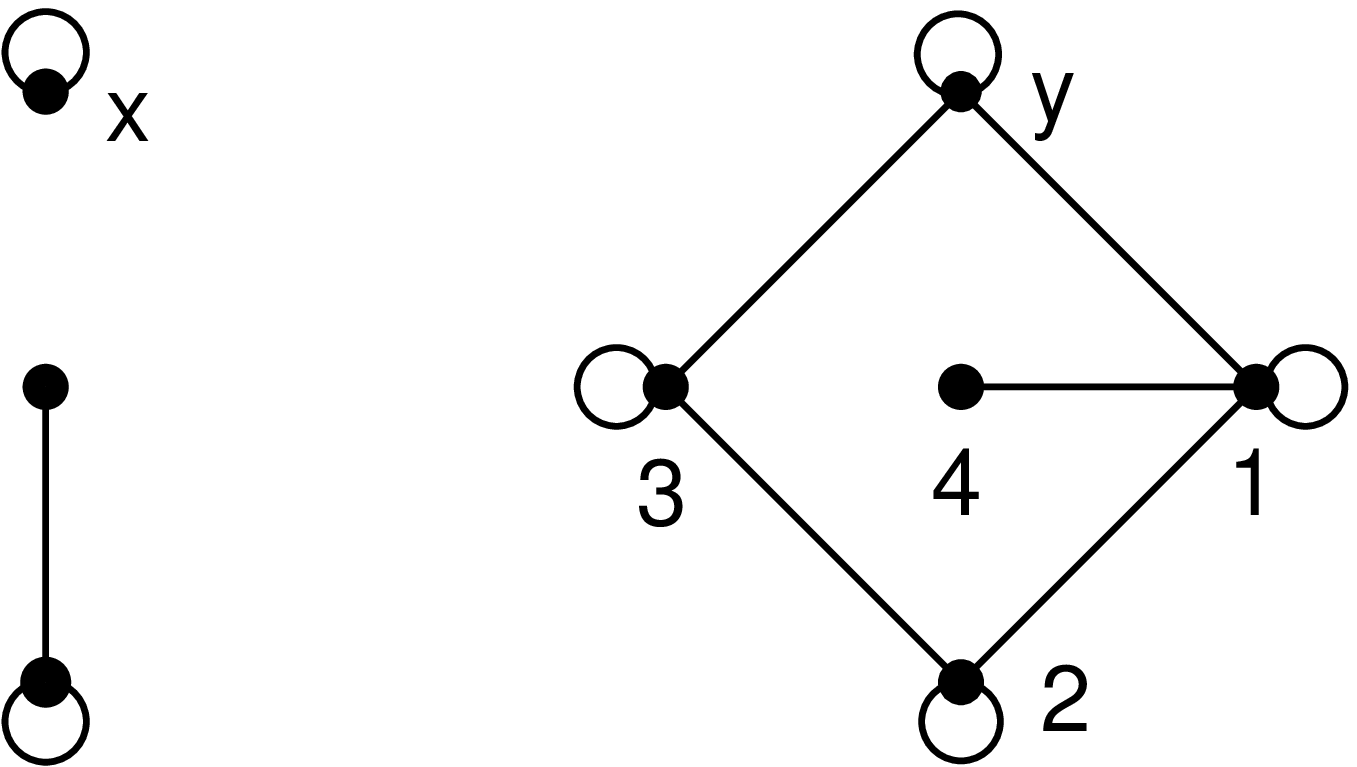, height=1.2 in, width = 2 in}

{Figure 3: The graphs $(G,x)$ and $(H,y)$.}
\end{center}

Each element of $\Hom_*(G,H)$ consists of certain functions $\eta$ from the vertex set of $G$ to nonempty subsets of the vertex set of $H$.  In particular, the pointed vertex $x \in V(G)$ must be sent to $\{y\} \subset V(H)$ for every $\eta$.  The realization of $\Hom_*(G,H)$ is the barycentric subdivision of the complex pictured in Figure 4, where the atoms of the poset are labeled with the images of the nonpointed vertices of $G$.  In this case we see that $\Hom_*(G,H) \simeq {\mathbb S}^1$.

\begin{center}
\epsfig{file=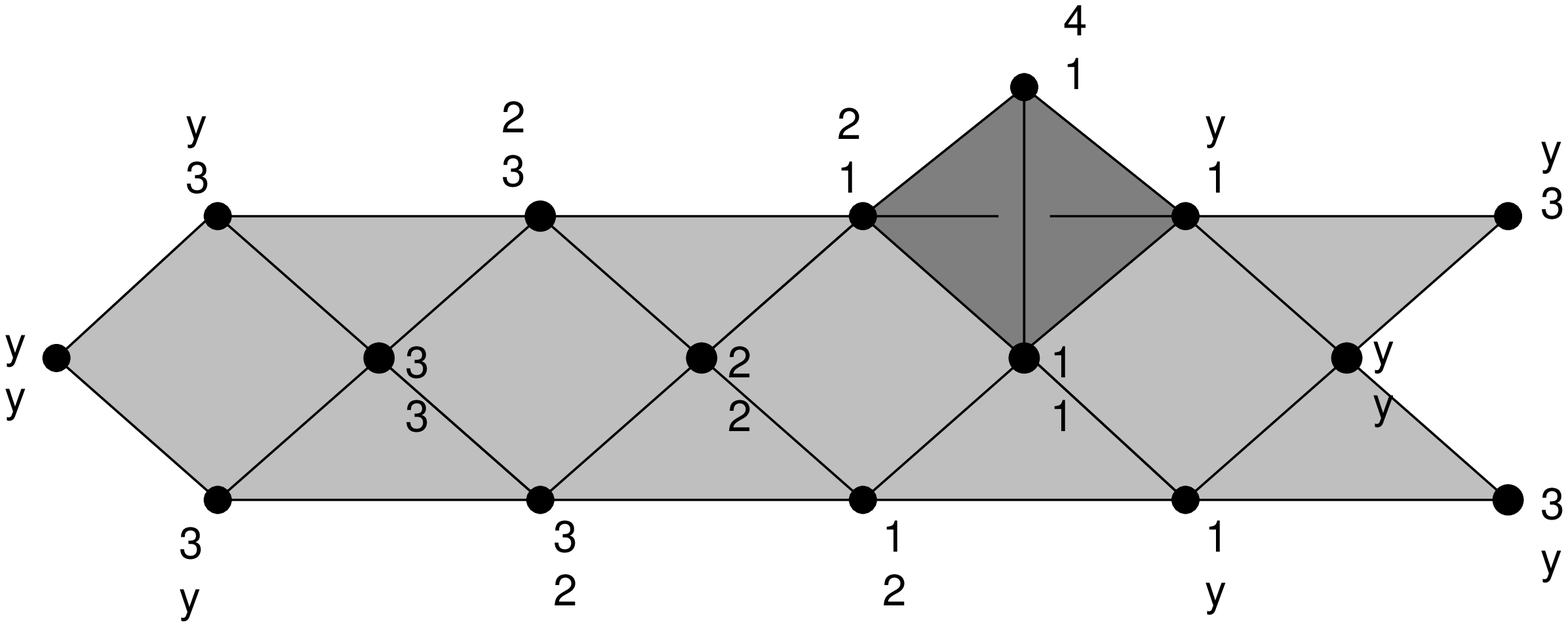, height=2 in, width = 4.3 in}

{Figure 4: The realization of $\Hom_*(G,H)$, up to barycentric subdivision.}
\end{center}

\end{example}

Note that $\Hom_*(G,H)$ is itself a pointed poset, with a distinguished element given by the map that sends all
vertices of $G$ to the vertex $y \in V(H)$.  One can also check that if $G$ is a connected graph with at least one edge, we have isomorphisms of posets 
\[\Hom_*(G_*,H_*) = \Hom(G, H_*) = \Hom(G,H) \coprod \Hom(G,*) = \big(\Hom(G,H)\big)_*,\]
where the last poset is obtained by adding a disjoint element to $\Hom(G,H)$.

For pointed graphs $A=(A,x)$, $B=(B,y)$, and  $C=(C,z)$, we have seen that
the exponential graph construction provides the adjunction
${\mathcal G}_*(A \wedge B, C) = {\mathcal G}_*(A, C^B)$, an
isomorphism of sets.  As is the case in the unpointed context (see
\cite{DocHom} and also \cite{Kfold}), this extends to a homotopy
equivalence of the analogous $\Hom_*$ complexes.

\begin{prop} \label{pointadjoint}
Let $A = (A,x),B = (B,y)$, and $C = (C,z) $ be pointed graphs.  The complex $\Hom_*(A \wedge B, C)$ can be included in $\Hom_*(A, C^B)$ so that
$\Hom_*(A \wedge B, C)$ is a strong deformation retract of
$\Hom_*(A, C^B)$. In particular, we have $\Hom_*(A \wedge B, C)
\simeq \Hom_*(A, C^B)$.
\end{prop}

\begin{proof}
We follow the proof of the analogous statement in the unpointed
context (see \cite{DocHom}).  For convenience we let $P = \Hom_*(A \wedge B, C)$
and $Q = \Hom_*(A, C^B)$ denote the respective posets.  We define a
map of posets $j: P \rightarrow Q$ according to
\[j(\alpha)(a) = \big\{f:\big(V(B),y\big) \rightarrow \big(V(C),z\big)|f(b) \in \alpha(a,b) ~\forall b \in B \big\},\]
\noindent
for every $\alpha \in P$ and $a \in A$.  Note that $z \in \alpha(x,b)$ for all $b \in V(B)$ and hence the
constant function $f_z:V(B) \rightarrow V(C)$, given by $b \mapsto
z$ for all $b$, is an element of $j(\alpha)(x)$. Also, if $(a,
a^\prime) \in E(A)$ then we have $(f,f^\prime) \in E(C^B)$ for all
$f \in j(\alpha)(a), f^\prime \in j(\alpha)(a^\prime)$. Hence
$j(\alpha)$ is indeed an element of $\Hom_*(A, C^B)$.

To see that $j$ is injective, suppose $\alpha \neq
\alpha^\prime \in \Hom_*(A \wedge B, C)$ with $\alpha(a,b) \neq
\alpha^\prime(a,b)$.  Then we have $\big\{f(b)| f \in j(\alpha)(a)\big\}
\neq \big\{g(b)| g \in j(\alpha^\prime)(a)\big\}$, so that indeed $j(\alpha) \neq
j(\alpha^\prime)$.

Next we define a closure operator $c: Q \rightarrow Q$.  If
$\gamma:V(A) \rightarrow 2^{V(C^B)} \backslash \{\emptyset \}$ is an
element of $Q$, we define $c(\gamma) \in \Hom_*(A, C^B)$ as
follows.  Fix some $a \in V(A)$ and for every $b \in V(B)$ let
$C_{ab}^\gamma = \big\{f(b) \in V(C)| f \in \gamma(a)\big\}$; define $c(\gamma)(a)$ to be the collection of functions $g:V(B) \rightarrow
V(C)$ where $g(b)$ varies over all $x \in C_{ab}^\gamma$.  One can verify that $c(\gamma) \in \Hom_*(A, C^B)$ and also that $c(p)
\geq p$ and $(c \circ c)(p) = c(p)$ for all $p \in P$.

Next we claim that $c(Q) \subseteq j(P)$.  To see this, suppose $\gamma \in Q = \Hom_*(A, C^B)$ so that $c(\gamma) \in c(Q)$.  We define $\alpha: V(A \wedge B)
\rightarrow 2^{V(C)} \backslash \{\emptyset \}$ by $\alpha(a,b) =
C_{ab}^\gamma$, where $C_{ab}^\gamma \subseteq V(C)$ is as
above. One can verify that $\alpha \in \Hom_*(A \wedge B, C)$.

Finally $j(P) \subseteq c(Q)$ since $j(P) \subseteq Q$ and $c(j(P)) =
j(P)$.  We conclude $j(P) = c(Q)$, so that $P \cong j(P)$ is the
image of a closure operator on $Q$.  Closure operators of posets induce strong deformation retracts of their order complexes (see for instance \cite{Bjo95}) and so the result follows.
\end{proof}

We let ${\bf 1}_*$ denote the graph consisting of a pair of
disjoint, looped vertices, and note that $G \wedge {\bf 1}_* = G$
for every pointed graph $G = (G,x)$.  The above proposition gives us
$\Hom_*(G,H) = \Hom_*({\bf 1}_* \wedge G, H) \simeq \Hom_*({\bf
1}_*, H^G)$, for pointed graphs $G = (G,x)$ and $H = (H,y)$. The
last of these posets is simply the face poset of the \textit{clique complex} of $H^G$, denoted $\Delta(H^G)$, which is by definition the simplicial complex whose faces are given by complete subgraphs on the looped vertices of $H^G$.  Recall that the
looped vertices in $H^G$ are the (pointed) graph homomorphisms
$(G,x) \rightarrow (H,y)$. Hence, for pointed graphs $G$ and $H$,
the complex $\Hom_*(G,H)$ can be realized up to homotopy type as the
clique complex of the subgraph of $H^G$ induced by the (pointed)
graph homomorphisms.

With this observation we obtain the following characterization of
$\times$-homotopy.

\begin{lemma} \label{pointcomponent}
Suppose $G=(G,x)$ and $H=(H,y)$ are pointed graphs, and $f,g:G
\rightarrow H$ are pointed graph maps. Then $f$ and $g$ are
$\times$-homotopic (as pointed maps) if and only if they are in the
same path-connected component of $\Hom_*(G,H)$.
\end{lemma}

\begin{proof}
A $\times$-homotopy $F:I_{n *} \rightarrow H^G$ from $f$ to $g$ is a
path in the 1-skeleton of $\Delta(H^G) \simeq \Hom_*(G,H)$.
Conversely, a (topological) path $I \rightarrow |\Hom_*(G,H)|$ can be approximated
as a simplicial map from some finite subdivision of $I$ into
$\Hom_*(G,H) \simeq \Delta(H^G)$.
\end{proof}

For a (not necessarily pointed) graph map $f:G \rightarrow H$, and a
fixed graph $T$, the $\Hom(T,?)$ and $\Hom(?,T)$ functors provide
maps $f_T$ and $f^T$ in the category of topological spaces. In Theorem 5.1 of \cite{DocHom} it is shown that $\times$-homotopy of graph maps is
characterized by the homotopy properties of these induced maps. For
our purposes, we will only need the following implication.

\begin{lemma} \label{homotopic}
Let $f,g:G \rightarrow H$ be maps of graphs.  Then $f$ and $g$ are
$\times$-homotopic if and only if the induced maps of posets $f_T, g_T:
\Hom(T,G) \rightarrow \Hom(T,H)$ are homotopic for every graph $T$.
\end{lemma}

The pointed $\Hom_*$ complexes also interact well with a graph
operation known as \textit{folding}.  We review this construction
next.

\begin{defn}
Let $u$ and $v$ be vertices of a pointed graph $G = (G,x)$ with $v
\neq x$ such that $N(v) \subseteq N(u)$.  Then we have a (pointed)
map $f:G \rightarrow G \backslash v$ given by $f(y) = y$, $y \neq
v$, and $f(v) = u$. We call the map $f$ a \textit{folding} of $G$ at
the vertex $v$.  The inclusion $i:G \backslash v \rightarrow G$
is called an \textit{unfolding}.
\end{defn}

\begin{prop} \label{pointfolds}
Suppose $G=(G,x)$ and $H=(H,y)$ are pointed graphs, and $u$ and $v$
are vertices of $G$ with $v \neq x$ such that $N(v) \subseteq N(u)$.
Then $i^H:\Hom_*(G,H) \rightarrow \Hom_*(G \backslash v, H)$ and
$f_H: \Hom_*(H,G) \rightarrow \Hom_*(H,G\backslash v)$ are both strong
deformation retracts, where $i^H$ and $f_H$ are the poset maps
induced by the graph unfolding and folding maps.
\end{prop}

\begin{proof}
We mimic the proof given in \cite{Kfold} of the analogous statement
in the unpointed setting.  For the first deformation retract, we
identify $\Hom_*(G \backslash v, H)$ with the subposet of
$\Hom_*(G,H)$ consisting of all $\alpha$ such that $\alpha(v) =
\alpha(u)$. We define $X$ to be the subposet of $\Hom_*(G,H)$ given by
$X = \{\alpha \in \Hom_*(G,H): \alpha(u) \subseteq \alpha(v)\}$.
Next, we define poset maps $\varphi:\Hom_*(G,H) \rightarrow X$ and
$\psi:X \rightarrow \Hom_*(G \backslash v, H)$ according to:

\begin{displaymath}
\varphi(\alpha)(w) = \left\{ \begin{array}{ll}
\alpha(u) \cup \alpha(v) & \textrm{if $w = v$}\\
\alpha(w) & \textrm{otherwise}
\end{array} \right.
\end{displaymath}

\begin{displaymath}
\psi(\alpha)(w) = \left\{ \begin{array}{ll}
\alpha(u) & \textrm{if $w = v$}\\
\alpha(w) & \textrm{otherwise.}
\end{array} \right.
\end{displaymath}

We see that both $\psi$ and $\varphi$ are closure maps, and that
$i^H = \psi \varphi$.  As above, closure maps of posets induce strong deformation retracts of their order complexes, and since $\im(i^H) = \Hom_*(G \backslash v, H)$ we obtain the result for $i^H$.

For the other statement, we define $Y$ to be the subposet of
$\Hom_*(H,G)$ given by $Y = \{\beta \in \Hom_*(H,G):\beta(w) \cap \{u,v\} \neq \{v\}$ for all $w \in V(H)\}$.  Define
poset maps $\rho:\Hom_*(H,G) \rightarrow Y$ and $\sigma:Y
\rightarrow \Hom_*(H,G \backslash v)$ according to:

\begin{displaymath}
\rho(\beta)(w) = \left\{ \begin{array}{ll}
\beta(w) \cup \{u\} & \textrm{if $v \in \beta(w)$}\\
\beta(w) & \textrm{otherwise}
\end{array} \right.
\end{displaymath}

\begin{displaymath}
\sigma(\beta)(w) = \left\{ \begin{array}{ll}
\beta(w) \backslash \{v\} & \textrm{if $v \in \beta(w)$}\\
\beta(w) & \textrm{otherwise.}
\end{array} \right.
\end{displaymath}

Once again, we see that both $\rho$ and $\sigma$ are closure maps,
with $\sigma \rho = f_H$.  Since $\im(f_H) = \Hom_*(H, G \backslash
v)$, the result follows.
\end{proof}

\begin{rem}
In \cite{Kfold}, Kozlov has shown that a closure map $c:P \rightarrow P$ induces a collapsing of the order complex of $P$ onto the order complex of $c(P)$.  Accordingly, one can strengthen the conclusions of Propositions \ref{pointadjoint} and \ref{pointfolds} to obtain a \textit{simple} homotopy equivalence between the relevant spaces.  We have only stated them in the form sufficient for our purposes.
\end{rem}

For each $n \geq 0$, the graph $I_n$ is pointed by the vertex $0$.
We can use Lemma \ref{pointfolds} to show that all $\Hom_*$
complexes involving the pointed graph $I_n$ are contractible.

\begin{lemma}\label{contractible}
For every pointed graph $G = (G,x)$ and for every integer $n$, the
complex $\Hom_*(I_n, G)$ is contractible.
\end{lemma}

\begin{proof}
If $n = 0$, then $I_0$ is a single looped vertex, and
$\Hom_*(I_0,G)$ is a vertex.  For $n > 0$, we have $N(n) \subset
N(n-1)$, and hence the unfolding map $i:I_{n-1} \rightarrow I_n$
induces a homotopy equivalence $f^G:\Hom_*(I_n,G) \rightarrow
\Hom_*(I_{n-1}, G)$.  The claim follows by induction.
\end{proof}

\section{$T$-homotopy groups and the main result}

We have seen that the path components of $\Hom_*(G,H)$ characterize
the $\times$-homotopy groups of maps from $G$ to $H$.  To relate the
higher homotopy groups to graph theoretical constructions we will
want to work with the `path graph' of a given pointed graph $G$,
which we define next.  For this and subsequent constructions we will make use of the notion of a \textit{colimit of a diagram} (of graphs or topological spaces).  We refer to \cite{Mac98} for a thorough discussion of this concept, but point out that in our context all such colimits are obtained from sequences of inclusions and hence can be thought of simply as unions.  For our first such construction, recall that the graph $I_n$ is pointed by the
vertex $0$.

\begin{defn}
For a pointed graph $G = (G,x)$, we define $G^I$ to be the pointed
graph obtained as the colimit (union) of the diagram,
\begin{center}

$\xymatrix{
G^{I_0} \ar[drr] \ar[r] & \dots \ar[r] & G^{I_{n}} \ar[d] \ar[r] & G^{I_{n+1}} \ar[dl] \ar[r] & \dots \\
& & G^I & & &} $

\end{center}

\noindent

where the maps $j_n: G^{I_n} \rightarrow G^{I_{n+1}}$ are induced by
the maps $I_{n+1} \rightarrow I_n$ given by $i \mapsto i$ ($i \neq
n+1$), and $n+1 \mapsto n$ (see Figure 5).

\end{defn}

\begin{center}

\epsfig{file=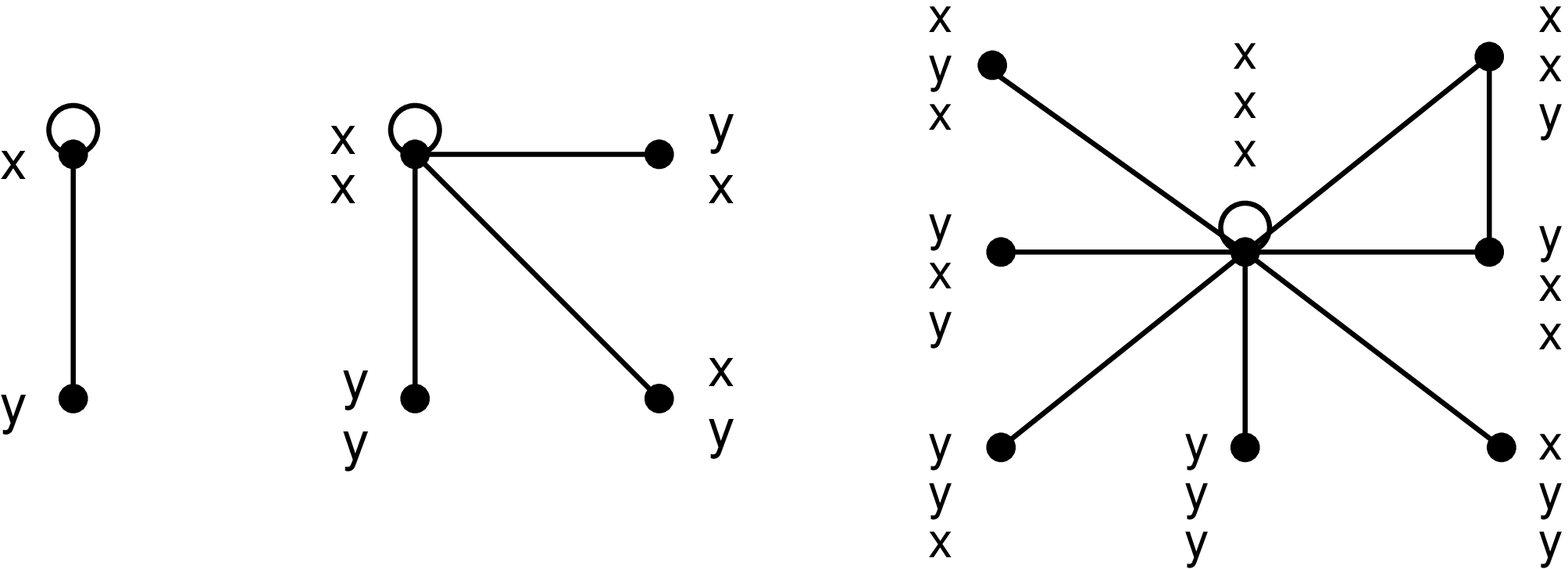, height=1.5 in, width = 3.8 in}

{Figure 5: The graphs $G$, $G^{I_1}$ (with the images of $0$, $1$),
and $G^{I_2}$ (with the images of $0$, $1$, $2$).}
\end{center}

Note that a vertex of $G^I$ is a (set) map $f:{\mathbb N}
\rightarrow V(G)$ from the nonnegative integers into the vertices of
$G$ with $f(0) = x$ which is \textit{eventually constant}; there exists some integer $N_f$ such that $f(i) = f(j)$ for all $i, j
\geq N_f$.  Adjacency is given by $f \sim g$ if $f(i) \sim g(j)$,
for all $i,j$ with $|i-j| \leq 1$. We think of $G^I$ as the graph
that parameterizes the collection of paths in $G = (G,x)$ which begin at the vertex $x \in G$.  The looped vertices of $G^I$ are those paths which involve only looped vertices of $G$.  We have the \textit{endpoint map} $\varphi:G^I \rightarrow G$ given by
$\varphi(f) = f(N_f)$.

\begin{defn}
For a pointed graph $G = (G,x)$, we define the \textit{loop space
graph} $\Omega G$ to be the (pointed) subgraph of $G^I$ induced by
elements that are eventually constant on the vertex $x \in G$.
\end{defn}

Hence $\Omega G$ is the graph whose vertices are given by closed paths in $G$ that start at $x$ and eventually end
(and stabilize) at $x$.  The looped vertices of $\Omega G$ are those closed paths that involve only looped vertices of $G$. The path and loop space graph functors commute with
exponentials of finite graphs, as described by the following observation.

\begin{lemma}\label{commute}
If $G=(G,x)$ and $T=(T,y)$ are pointed graphs, with $T$ finite, then
we have graph isomorphisms $(G^I)^T = (G^T)^I$ and $(\Omega G)^T = \Omega(G^T)$.
\end{lemma}

\begin{proof}
Both isomorphisms follow from identical arguments, and so we prove
only the second of these claims.  Define a map $\alpha:(\Omega G)^T
\rightarrow \Omega(G^T)$ by $\alpha(f)(i)(t) = f(t)(i)$, for $f \in
(\Omega G)^T$.  To show that $\alpha(f) \in \Omega(G^T)$, pick an
integer $j$ such that each element of \{$f(t)\}_{t \in V(T)}$
stabilizes at $j$ (this is possible since $V(T)$ is finite).  So we
have $f(t)(k) = f(t)(k^\prime) = x$ for all $t \in V(T)$ and $k,
k^\prime \geq j$. Hence $\alpha(f)(k) = \alpha(f)(k^\prime)$ so that
$\alpha(f)$ stabilizes at $j$.  It is easy to check that $\alpha$ is
a pointed graph map. To show that it is an isomorphism, define a map
$\beta:\Omega(G^T) \rightarrow (\Omega G)^T$ by $\beta(g)(t)(i) =
g(i)(t)$.   Suppose $g \in \Omega(G^T)$ stabilizes at the integer
$j$.  Then we have $g(k) = g(k^\prime)$ as elements of $G^T$, for
all $k, k^\prime \geq j$.  Hence for every $t$, we have that
$\beta(g)(t)(k) = \beta(g)(t)(k^\prime)$ so that in fact $\beta$
maps to $(\Omega G)^T$.  Again one can check that $\beta$ is a
pointed graph map and the inverse to $\alpha$.  The result follows.
\end{proof}

Note that $\Omega G$ is pointed by the
closed path that is constant on the vertex $x \in G$.  
More generally, there is a natural group structure on the set $\pi_0 \big(\Hom(T, \Omega H)\big)$ for any graph $H$ and finite graph $T$.  To describe this structure, we let $G := H^T$ and appeal to Lemma \ref{commute} and Proposition \ref{pointadjoint} to obtain the string of bijections $\pi_0\big(\Hom_*(T,\Omega H)\big) \cong \pi_0\big(\Hom_*({\bf 1}_*, (\Omega H)^T)\big) \cong \pi_0\big(\Hom_*({\bf 1}_*, \Omega(H^T))\big) \cong [{\bf 1}_*, \Omega G]_{\times}$.  Hence for our purposes it is sufficient to describe a group structure on the set $[{\bf 1}_*, \Omega G]_{\times}$, whose elements we can identify with connected components of $\Omega G$.

We define a multiplication on the components of $\Omega$ in the following way.  Given a pair of elements $f$ and $g$ in $\Omega G$, pick a number $N_g$ such that $g(n) = x$ for all $n \geq N_g$.  Define $[f] \cdot [g] := [(f \cdot_{N_g} g)]$, where 
\begin{displaymath}
(f \cdot_{N_g} g)(i) = \left\{ \begin{array}{ll}
g(i) & \textrm{if $i < N_g$}\\
f(i-N_g) & \textrm{otherwise.}
\end{array} \right.
\end{displaymath}
The identity element is given by $[c_x]$, where $c_x$ is the path that is constant on the vertex $x \in G$.  The inverse of $[f]$ is given by (the equivalence class of) the path $f$ traversed in the opposite direction, so that
\begin{displaymath}
[f]^{-1}(i) = \left\{ \begin{array}{ll}
f(N_f - i) & \textrm{if $i \leq N_f$}\\
x & \textrm{otherwise.}
\end{array} \right.
\end{displaymath}

First we show that our construction does not depend on the choice of $N_g$.  For this, suppose we pick $N_g^\prime > N_g$ in our construction of $[f] \cdot [g]$; we need to show that $[f] \cdot_{N^\prime} [g]$ and $[f] \cdot_{N} [g]$ are in the same component of $\Omega G$.  By induction it is enough to assume that $N_g^\prime = N_g + 1$, and in this case we build a path $(f \cdot_{N^\prime} g = f \cdot_{N+1} g := h_0, h_1, \dots, h_{k-1}, h_k := f \cdot_{N} g)$ in $\Omega G$ according to 
\begin{displaymath}
h_j(i) = \left\{ \begin{array}{ll}
g(i) & \textrm{if $i < N_g$}\\
f(i-N_g) & \textrm{if $N_g \leq i < N_g + j$}\\
f(i-N_g-1) & \textrm{if $i \geq N_g + j$}
\end{array} \right.
\end{displaymath}

This path is perhaps best understood with the help of the following diagram, where each row represents an element $h_j$, and each column represents the image of $i$ under $h_j$.

\begin{center}
\begin{tabular}{ccccccccc}
$N_g$ & $N_g + 1$ & $N_g+2$ & $N_g+3$ & & $\cdots$ & $N_g + N_f -1$ & $N_g + N_f$ & $N_g + N_f + 1$ \\
\hline
$x$ & $x$ & $f(1)$ & $f(2)$ & $f(3)$ & $\cdots$ & $f(N_f - 2)$ & $f(N_f - 1)$ & $f(N_f) = x$ \\
$x$ & $f(1)$ & $f(1)$ & $f(2)$ & $f(3)$ & $\cdots$ & $f(N_f - 2)$ & $f(N_f - 1)$ & $x$ \\
$x$ & $f(1)$ & $f(2)$ & $f(2)$ & $f(3)$ & $\cdots$ & $f(N_f - 2)$ & $f(N_f - 1)$ & $x$ \\
$x$ & $f(1)$ & $f(2)$ & $f(3)$ & $f(3)$ & $\cdots$ & $f(N_f - 2)$ & $f(N_f - 1)$ & $x$ \\
$\cdots$ & $\cdots$ & $\cdots$ & $\cdots$ & $\cdots$ & $\cdots$ & $\cdots$ & $\cdots$ & $\cdots$\\
$x$ & $f(1)$ & $f(2)$ & $f(3)$ & $f(4)$ & $\cdots$ & $f(N_f - 1)$ & $f(N_f - 1)$ & $x$ \\
$x$ & $x$ & $f(1)$ & $f(2)$ & $f(3)$ & $\cdots$ & $f(N_f - 1)$ & $f(N_f)=x$ & $x$
\end{tabular}
\end{center}

We conclude that $f \cdot_{N} g$ and $f \cdot_{N^\prime} g$ are in the same component of $\Omega G$, and hence our product does not depend on the choice of $N_G$ (up to $\times$-homotopy).  Similarly, if $f \simeq_{\times} f^\prime$ and $g \simeq_{\times} g^\prime$ in $[{\bf 1}_*,\Omega G]_{\times}$ (so that $f$ and $f^\prime$, resp. $g$ and $g^\prime$, are in the same component of $\Omega G$), one can check that $f \cdot g$ and $f^\prime \cdot g^\prime$ are in the same component of $\Omega G$, and hence the product is well defined on $\times$-homotopy classes of $[{\bf 1}_*,\Omega G]_{\times}$.  It is clear that $[c_x] \cdot [f] = [f] \cdot [c_x] = [f]$, and it is not hard to see that $[f] \cdot [f]^{-1} = [c_x] = [f]^{-1} \cdot [f]$.  

As an example of the latter claim, consider the graph $(H,y)$ in Figure 3, and let $f = (y,1,2,3,y,y,\dots)$ be an element of $\Omega H$.  Taking $N_f = 4$, we obtain $(y,3,2,1,y,y,\dots)$ as an element of the equivalence class $[f]^{-1}$, and $(y,1,2,3,y,y,3,2,1,y,y,\cdots)$ as an element of $[f]^{-1} \cdot [f]$.  We get a homotopy from this product path to $c_x$ according to 
\begin{center}
\begin{tabular}{ccccccccccc}
0 & 1 & 2 & 3 & 4 & 5 & 6 & 7 & 8 & 9 & $\cdots$\\
\hline
$y$ & 1 & 2 & 3 & $y$ & $y$ & 3 & 2 & 1 & $y$ & $y$\\
$y$ & 1 & 2 & 3 & 3   &  3  & 3 & 2 & 1 & $y$ & $y$\\
$y$ & 1 & 2 & 2 & 2   &  2  & 2 & 2 & 1 & $y$ & $y$\\
$y$ & 1 & 1 & 1 & 1   &  1  & 1 & 1 & 1 & $y$ & $y$\\
$y$ & $y$ & $y$ & $y$ & $y$ & $y$ & $y$ & $y$ & $y$ & $y$ & $y$
\end{tabular}
\end{center}
 
As one might expect, the group structure described above coincides with the fundamental group of the complex $\Hom_*(T,G)$.

\begin{lemma}\label{fundgroup}
Let $G = (G,x)$ and $T = (T,y)$ be pointed graphs, with $T$ finite.  Then we have an isomorphism of groups
\[\pi_1\big(\Hom_*(T,G)\big) \cong \pi_0\big(\Hom_*(T, \Omega G)\big),\]
\noindent
where the basepoint of $\Hom_*(T,G)$ is taken to be the element which sends each vertex of $T$ to $\{x\}$, the set containing the pointed vertex of $G$.
\end{lemma}

\begin{proof}
Each element $\alpha \in \pi_1\big(\Hom_*(T,G)\big) \cong \pi_1\big(\Hom_*({\bf 1}_*, G^T)\big)$ can be represented up to homotopy by a simplicial map from a subdivision of the circle (say with $n$ vertices) into the clique complex of the graph $G^T$.  This map can be thought of as a pointed graph map $\tilde \alpha :C^\prime_n \rightarrow G^T$, where $C^\prime_n$ is a reflexive (loops on all vertices) cycle of length $n$ (the graph $I_n$ with the endpoints identified) and where the basepoint of $G^T$ is the vertex set map which sends all vertices of $T$ to $x \in G$.  Hence $\alpha$ can be considered to be an element of $\Omega(G^T)$.  Similarly, if $\alpha$ and $\beta$ are homotopic as pointed maps ${\mathbb S}^1 \rightarrow G^T$, then for sufficiently large $n$ these can be represented as graph maps $\tilde \alpha, \tilde \beta: C^\prime_n \rightarrow G^T$.  The maps $\alpha$ and $\beta$ are recovered (up to homotopy) by the induced $\tilde \alpha_{{\bf 1}_*}, \tilde \beta_{{\bf 1}_*}: \Hom_*({\bf 1}_*, C^\prime_n) \rightarrow \Hom_*({\bf 1}_*, G^T)$.  These maps are homotopic by assumption, and hence by Lemma \ref{homotopic} we have that $\tilde \alpha$ and $\tilde \beta$ are $\times$-homotopic as graph maps.  By definition this means that we have a path along looped vertices in $(G^T)^{C^\prime_n}$ between (the looped vertices) $\tilde \alpha$ and $\tilde \beta$.  But $(G^T)^{C^\prime} \subseteq \Omega(G^T)$ and hence this provides a path in $\Omega(G^T)$.   This in turn implies that $\tilde \alpha$ and $\tilde \beta$ are in the same component of $\Omega(G^T)$ and hence the same as elements of $\pi_0\big(\Hom_*({\bf 1}_*, \Omega(G^T))\big) \cong \pi_0\big(\Hom_*(T,\Omega G)\big)$. 

This assignment is a homomorphism of groups since concatenation of loops corresponds to the group structure on $\Omega (G^T)$ described above.  Hence we obtain a map $\varphi:\pi_1\big(\Hom_*({\bf 1}_*, G^T)\big) \rightarrow \pi_0\big(\Hom_*(T, \Omega G)\big)$.  As above, any element of $\pi_0 \big(\Hom_*(T, \Omega G)\big)$ can be represented by a graph map from some finite reflexive cycle $C^\prime_n$ into $G^T$, and hence $\varphi$ is surjective.  For injectivity of $\varphi$, suppose $\alpha$ and $\beta$ are elements in the same $\times$-homotopy class of $[{\bf 1}_*, \Omega(G^T)]_{\times}$.  We can choose $n$ large enough so that $\alpha$ and $\beta$ are realized as $\times$-homotopic elements of $[C^\prime_n, G^T]_{\times}$.  Applying $\Hom_*({\bf 1}_*,?)$ to these maps then gives representives in $\pi_1\big(\Hom_*({\bf 1}_*, G^T)\big)$ which are homotopic by Lemma \ref{homotopic}. 
\end{proof}

We next turn to a consideration of the higher homotopy groups of $\Hom_*$ complexes.  For this we will consider iterations of the loop space construction, and will use $\Omega^n (G)$ to denote the graph $\Omega \big(\Omega \big(\dots \big(\Omega(G) \big)\dots \big) \big)$ ($n$ times).  We will also need the following observation.  

\begin{lemma} \label{pathcontract}
For a pointed graph $G$ and a finite pointed graph $T$,
$\Hom_*(T,G^{I})$ is contractible.
\end{lemma}

\begin{proof}
We prove the claim for $T = {\bf 1}_*$, from which the result
follows from the homotopy equivalence $\Hom_*(T,G^{I}) \simeq
\Hom_*({\bf 1}_*,(G^{I})^T) = \Hom_*({\bf 1}_*, (G^T)^I)$.

We first show that $\Hom_*({\bf 1}_*,G^{I_n})$ is contractible for
any integer $n$.  By Proposition \ref{pointadjoint} we have that
$\Hom_*({\bf 1}_*,G^{I_n}) \simeq \Hom_*({\bf 1}_* \wedge I_n, G) =
\Hom_*(I_n, G)$ (where $I_n$ is pointed by the vertex $0$).  The
latter space is contractible by Lemma \ref{contractible}.

Finally we prove that $\Hom_*({\bf 1}_*, G^I)$ is contractible.  We
have seen that $\Hom_*({\bf 1}_*, X) = \Hom({\bf 1}, X)$ is the
clique complex on the looped vertices of the graph $X$. Hence, as a
functor, $\Hom_*({\bf 1}_*, ?)$ preserves colimits. By definition,
$G^I$ is the colimit of the sequence of maps $\dots \rightarrow
G^{I_n} \rightarrow G^{I_{n+1}} \rightarrow \dots$, and hence
$\Hom_*({\bf 1}_*, G^I) = \colim \big(\Hom_*({\bf 1}_*,
G^{I_n})\big)$. We have seen that the sequence defining the last of
these spaces is composed of all contractible spaces. Hence the
colimit is contractible, and the result follows.
\end{proof}

Next we state and prove our main result, a long exact sequence in
the homotopy groups of the $\Hom_*$ complexes induced by the
endpoint map $\varphi:G^I \rightarrow G$.  Our main tool will be the
so-called Quillen fiber Lemma B (see \cite{Qui73}).  If $\psi:P
\rightarrow Q$ is a map of posets and $q \in Q$, we let
$\psi^{-1}(\leq q) := \{p \in P| \psi(p) \leq q \}$ denote the
\textit{Quillen fiber} of $q$.  

\begin{thm}[Quillen]
Let $\psi:P \rightarrow Q$ be a map of posets such that for all $q
\leq q^\prime$ the induced map $\psi^{-1}(\leq q) \rightarrow
\psi^{-1}(\leq q^\prime)$ is a homotopy equivalence.  Then for all $q \in Q$ and all $p \in \psi^{-1}(\leq q)$ there exists a connecting homomorphism $\delta:\pi_{i+1}(Q,q) \rightarrow \pi_i \big(\psi^{-1}(\leq q),p \big)$ that fits into the long exact sequence

\begin{center}
$\xymatrix{ \dots \ar[r] & \pi_{i+1}(Q,q) \ar[r]^\delta & \pi_i(\psi^{-1}(\leq q), p) \ar[r]^{\iota_*}  & \pi_i(P,p) \ar[r]^{\psi_*} & \pi_i(Q,q) \ar[r] &
\dots }$,
\end{center}

\noindent
where $\psi_*$ is the map induced by $\psi$ and $\iota_*$ is
induced by the inclusion $\iota:\big(\psi^{-1}(\leq q), p \big)
\rightarrow (P,p)$.
\end{thm}

\begin{thm} \label{sequence}
Let $G$ be a pointed graph, let $T$ be a finite pointed graph, and
let $\varphi_T: \Hom_*(T,G^I) \rightarrow \Hom_*(T,G)$ be the map
induced by the endpoint map.  Then for all $\gamma \in \Hom_*(T,G)$,
and all $\beta \in \varphi_T^{-1}(\leq \gamma)$ we have a connecting
homomorphism $\delta:\pi_{i+1}\big(\Hom_*(T,G),\gamma \big)
\rightarrow \pi_i \big(\varphi_T^{-1}(\leq \gamma), \beta \big)$ that fits
into the following long exact sequence.

\begin{center}
$\xymatrix{ \dots \ar[r] & \pi_{i+1}\big(\Hom_*(T,G),\gamma\big)
\ar[r]^\delta & \pi_i\big(\varphi_T^{-1}(\leq \gamma), \beta\big)
\ar[dl]_{\iota_*}  \\
& ~\pi_i\big(\Hom_*(T,G^I), \beta\big)
\ar[r]^{\varphi_*} & \pi_i\big(\Hom_*(T,G), \gamma \big) \ar[r] &
\dots }$
\end{center}

\noindent
Here $\varphi_*$ is the map induced by $\varphi_T$ and $\iota_*$ is
induced by the inclusion $\iota:(\varphi_T^{-1}(\leq \gamma), \beta)
\rightarrow (\Hom_*(T,G^I),\beta)$.
\end{thm}

\begin{proof}
We first prove the claim for the case $T = {\bf 1}_*$, and the map
$\varphi_{{\bf 1}_*}: \Hom({\bf 1}_*, G^I) \rightarrow \Hom({\bf
1_*}, G)$.

We use the Quillen fiber Lemma B applied to the poset map
$\varphi_{{\bf 1}_*}$.  Suppose $\gamma \leq \gamma^\prime \in
\Hom_*({\bf 1}_*, G)$.  As elements of $\Hom_*({\bf 1}_*, G)$, $\gamma$ and $\gamma^\prime$ can each be identified with a collection
of looped vertices of the graph $G$, each of which determines a
clique (complete subgraph) of $G$.  We will also use $\gamma$ and
$\gamma^\prime$ to denote these collections of vertices, and will
distinguish an element $v \in \gamma \subseteq \gamma^\prime
\subseteq V(G)$.  Note that $v$ is adjacent to all other elements of
$\gamma^\prime$ (see Figure 6).

\begin{center}
\epsfig{file=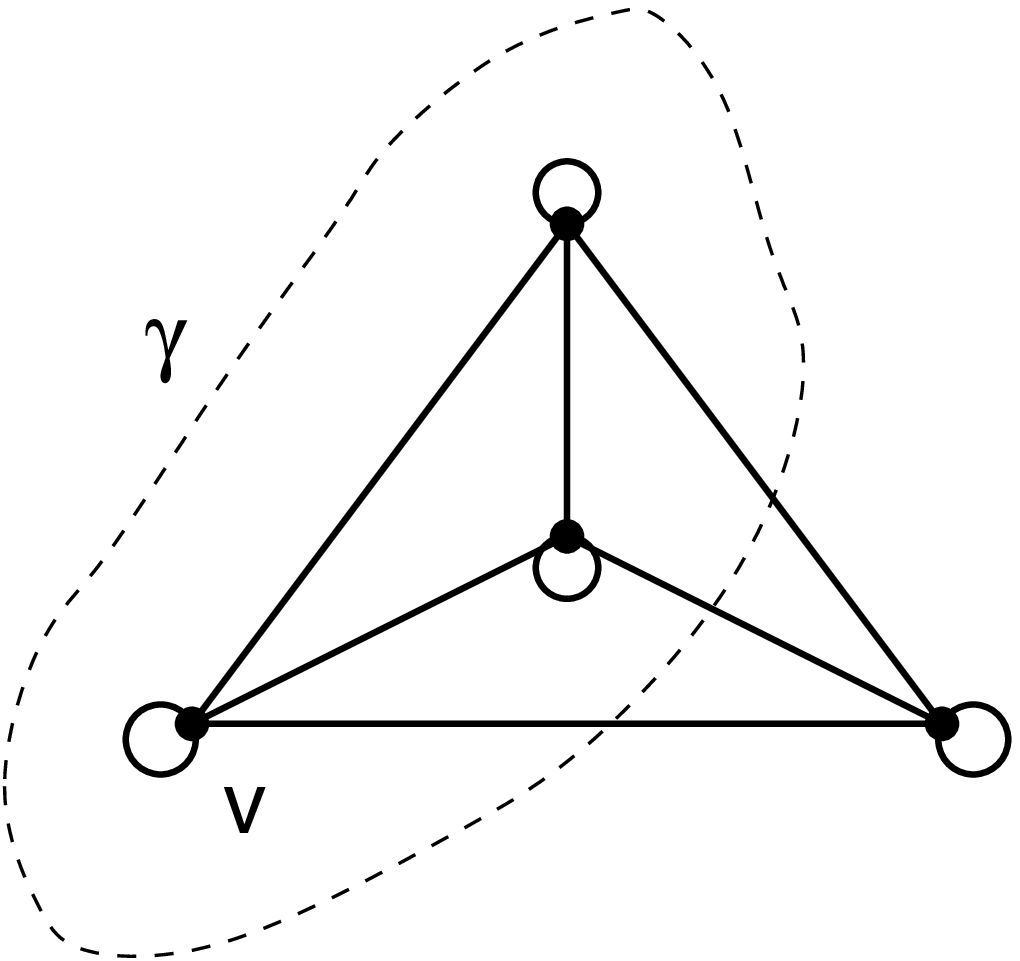, height=1.25 in, width = 1.25 in}

{Figure 6: $v \in \gamma \subseteq \gamma^\prime$.}
\end{center}

Next, we let $Y = \varphi^{-1}_{{\bf 1}_*}(\leq \gamma)$ and
$Y^\prime = \varphi^{-1}_{{\bf 1}_*}(\leq \gamma^\prime)$ denote the
respective Quillen fibers. Here $Y$ and $Y^\prime$ are both
subposets of $\Hom_*({\bf 1}_*, G^I)$.  We wish to show that the
induced map $k:Y \rightarrow Y^\prime$ is a homotopy equivalence,
from which the result would follow.  For this, we consider finite
approximations of these spaces, in the following sense.  We let
$\varphi_n: G^{I_n} \rightarrow G$ denote the endpoint map of the
finite path graph $G^{I_n}$, given by $\varphi_n(f) = f(n)$.  We let
$H_n \subseteq G^{I_n}$ denote the induced subgraph on the vertices
$\varphi_n^{-1}(\gamma)$, and similarly $H^\prime_n$ on
$\varphi_n^{-1}(\gamma^\prime)$.  We can think of $H_n$ (respectively
$H^\prime_n$) as the subgraph of $G^{I_n}$ induced by maps from
$I_n$ that end at some vertex of $\gamma$ (respectively
$\gamma^\prime$). We have the obvious inclusions $k_n:H_n
\rightarrow H^\prime_n$ and also the inclusions $i_n:H^\prime_n
\rightarrow H^\prime_{n+1}$ and $j_n:H_n \rightarrow H_{n+1}$ given
by $i_n(f)(n+1) = f(n)$ (and similarly for $j_n$).

We let $Y_n = \Hom({\bf 1}, H_n)$ and $Y^\prime_n = \Hom({\bf
1}, H^\prime_n)$ denote the respective posets. The $i_n$ and $j_n$
maps determine directed systems for which

\begin{center}
$Y = \varphi^{-1}(\leq \gamma) = \Hom({\bf 1}, \colim H_n ) = \colim Y_n$

$Y^\prime = \varphi^{-1}(\leq \gamma^\prime) = \Hom({\bf 1}, \colim
H^\prime_n) = \colim Y^\prime_n$.
\end{center}

\noindent
The poset map $k:Y \rightarrow Y^\prime$ is given by $\colim(k_n:H_n
\rightarrow H^\prime_n)$.  We also need the graph map
$h_n:H^\prime_n \rightarrow H_{n+1}$ given by
\begin{displaymath}
h_n(f)(i) = \left\{ \begin{array}{ll}
f(i) & \textrm{if $i \leq n$}\\
v & \textrm{if $i=n+1$}
\end{array} \right.
\end{displaymath}

These maps all fit into the following diagram of graphs.

\begin{center}
$\xymatrix{
\dots \ar[r] & H_n^\prime \ar[dr]^{h_n} \ar[r]^{i_n} & H_{n+1}^\prime \ar[dr]^{h_{n+1}} \ar[r]^{i_{n+1}} & H^\prime_{n+2} \ar[r] & \dots \\
\dots \ar[r] & H_n \ar[u]_{k_n} \ar[r]_{j_n} & H_{n+1}
\ar[u]_{k_{n+1}} \ar[r]_{j_{n+1}} & H_{n+2} \ar[u]_{k_{n+2}} \ar[r]
& \dots }$
\end{center}

We claim that this diagram commutes up to (graph) $\times$-homotopy.
In particular, we have $k_{n+1}h_n \simeq_\times i_n$ and
$h_{n+1}k_{n+1} \simeq_\times j_{n+1}$.  For the first homotopy,
define a map $A:H^\prime_n \times I_{1} \rightarrow H^\prime_{n+1}$
according to:
\begin{displaymath}
A(f,i)(j) = \left\{ \begin{array}{lll}
f(j) & \textrm{if $i=0$ and $j\leq n$}\\
f(n) & \textrm{if $i=0$ and $j=n+1$}\\
v & \textrm{if $i = 1$}
\end{array} \right.
\end{displaymath}
Recall that $v \in \alpha \subseteq V(G)$ is our distinguished
vertex.  Now, it is easy to see that $A(?,0) = i_n$ and $A(?,1) =
k_{n+1}h_n$.

\begin{center}
$\xymatrix{
H^\prime_n \ar[d]_{\iota_0} \ar[dr]^{i_n} \\
H^\prime_n \times I_{1} \ar[r]^A & H^\prime_{n+1} \\
H^\prime_n \ar[u]^{\iota_1} \ar [ur]_{k_{n+1}h_n} }$
\end{center}

To check that $A$ is a graph map, suppose $f$ and $f^\prime$ are adjacent vertices in $H^\prime_n$. We
need $A(f,0)$ and $A(f^\prime,1)$ to be adjacent in
$X^\prime_{n+1}$.  Note that $A(f,0)(n+1) = f(n) \in
\gamma^\prime$ and $A(f^\prime,1)(n+1) = v$.  The element $v$ is
adjacent to all elements in the clique $\gamma^\prime$, and hence
adjacent to $f(n)$. Also, $A(f,0)(n+1) = f(n)$ and
$A(f^\prime, 1)(n) = f^\prime(n)$, which are adjacent since
$f \sim f^\prime$ in $H^\prime_n$.  Finally, we have
$A(f,0)(n) = f(n)$ and $A(f,1)(n+1) = x$.  Once again,
these are adjacent since $f(n) \in \gamma^\prime$.

To check the homotopy $h_{n+1}k_{n+1} \simeq_{\times} j_{n+1}$, we similarly
define a map $B:H_n \times I_{1} \rightarrow H_{n+1}$ according to:
\begin{displaymath}
B(f,i)(j) = \left\{ \begin{array}{lll}
f(n) & \textrm{if $i=0$ and $j=n+1$}\\
v & \textrm{if $i=1$ and $j=n+1$}\\
f(j) & \textrm{if $j \leq n$}
\end{array} \right.
\end{displaymath}
Again one can check that $B$ is indeed a graph map and that
$B(?,0) = j_{n+1}$ and $B(?,1) = h_{n+1}k_{n+1}$.  We conclude that
the diagram under consideration commutes up to $\times$-homotopy,
and hence by Lemma \ref{homotopic} any diagram of posets induced by
a $\Hom(T,?)$ functor also commutes up to (topological) homotopy.

Next, to show that $k: Y \rightarrow Y^\prime$ is a homotopy
equivalence we show that $k$ induces an isomorphism on homotopy
groups.  Suppose $\rho, \sigma:{\mathbb S}^m \rightarrow Y$ are pointed maps
from the $m$-sphere into $Y$, and let $\rho^\prime = k \rho$ and
$\sigma^\prime = k \sigma$ be the induced maps ${\mathbb S}^m \rightarrow
Y^\prime$. Suppose that $\rho^\prime \simeq \sigma^\prime$ are
homotopic as maps into $Y^\prime$ via a homotopy $\Psi:{\mathbb S}^m \times I
\rightarrow Y^\prime$.  We claim that in fact $\rho \simeq \sigma$
are also homotopic, so that $k$ is injective on all homotopy groups.
To see this, pick $n$ big enough so that the image of $\Psi$ sits
inside the subcomplex $Y^\prime_n := \Hom({\bf 1}, H_n) \subseteq
Y^\prime$ (this is possible since ${\mathbb S}^m \times I$ is compact). Now,
the composition $h_{n_{\bf 1}} \Psi : {\mathbb S}^m \times I \rightarrow
Y_{n+1}$ is a homotopy from $h_{n_{\bf 1}} \rho^\prime = h_{n_{\bf
1}} k_{n_{\bf 1}} \rho$ to $h_{n_{\bf 1}} \sigma^\prime = h_{n_{\bf
1}} k_{n_{\bf 1}} \sigma$. But $h_{n_{\bf 1}} k_{n_{\bf 1}} \simeq
j_{n_{\bf 1}}$ and hence $\rho$ and $\sigma$ are homotopic as maps
into $Y_{n+1}$, as desired.

\begin{center}
$\xymatrix{ {\mathbb S}^m \ar[d] \ar[dr]_{\rho} & Y^\prime_n \ar[dr]^{h_{n_{\bf 1}}} \ar[r]^{i_{n_{\bf 1}}} & Y^\prime_{n+1} \\
{\mathbb S}^m \times I \ar[r]_{\Psi} & Y_n \ar[u]^{k_{n_{\bf 1}}} \ar[r]_{j_{n_{\bf 1}}} & Y_{n+1} \ar[u]_{k_{n_{\bf 1}}} \\
{\mathbb S}^m \ar[u] \ar[ur]_{\sigma} & & }$ 
\end{center}

Next, we claim that $k$ induces a surjection on each homotopy group.
To see this, suppose $\rho^\prime:{\mathbb S}^m \rightarrow Y^\prime$ is a
pointed map of the $m$-sphere into $Y^\prime$.  We wish to find a
map $\rho:{\mathbb S}^m \rightarrow Y$ such that $k \rho \simeq \rho^\prime$.
As above, choose $n$ large enough so that the image of the map
$\rho^\prime$ is contained in $Y^\prime_n$, and let $\rho =
h_{n_{\bf 1}} \rho^\prime: {\mathbb S}^m \rightarrow Y_{n+1}$.  Then we have
$k \rho = k_{(n + 1)_{\bf 1}} h_{n_{\bf 1}} \rho^\prime \simeq
i_{n_{\bf 1}} \rho^\prime \simeq \rho^\prime$, as desired. We
conclude that $k$ induces an isomorphism on each homotopy group, so
that $k$ is a homotopy equivalence.  Hence the conditions of the
Quillen Lemma $B$ are satisfied, and we get a long exact sequence

\begin{center}
$\xymatrix{ \dots \ar[r] & \pi_{i+1}\big(\Hom_*({\bf
1}_*,G),\gamma\big)
\ar[r] & \pi_i\big(\varphi_{{\bf 1}*}^{-1}(\leq \gamma), \beta\big) \ar[dl] \\
& \pi_i\big(\Hom_*({\bf 1}_*,G), \beta\big) \ar[r] &
\pi_i\big(\Hom_*({\bf 1}_*,G), \gamma\big) \ar[r] & \dots }$
\end{center}

It remains for us to prove the claim for general finite $T$.  For this we
note that $\Hom_*(T,G) \simeq \Hom_*({\bf 1}_*, G^T)$, and also
$\Hom_*(T, G^I) \simeq \Hom_*\big({\bf 1}_*, (G^I)^T\big) =
\Hom_*\big({\bf 1}_*, (G^T)^I\big)$.  Hence we take $G = G^T$ and
make the appropriate substitutions in the above sequence.
\end{proof}

\begin{cor}\label{Homgroups}
For pointed graphs $G$ and $T$, with $T$ finite, we have
$\pi_i\big(\Hom_*(T,G),\gamma \big) \simeq [T, \Omega^i(G)]_\times$.
\end{cor}

\begin{proof}
Suppose $T=(T,y)$ and $G=(G,x)$ are pointed graphs, with $T$ finite.
Let $\varphi:\Hom_*(T,G^I) \rightarrow \Hom_*(T,G)$ be the map of
posets induced by the (pointed) endpoint map $G^I \rightarrow G$. We
apply Theorem \ref{sequence} and in the long exact sequence choose
$\gamma \in \Hom_*(T,G)$ to be the basepoint (where $\gamma$ is
given by $\gamma(t) = x$ for all $t \in V(T)$).  Hence $\gamma$ is
an atom in the poset $\Hom_*(T,G)$, so that $\varphi_T^{-1}(\leq
\gamma) = \varphi_T^{-1}(\gamma) = \Hom_*(T, \Omega G)$.  We choose
$\beta$ to be the basepoint in $\Hom_*(T, \Omega G)$.  Hence our
sequence becomes

\begin{center}
$\xymatrix{ \dots \ar[r] & \pi_{i+1}\big(\Hom_*(T,G),\gamma \big)
\ar[r]^\delta  & \pi_i \big(\Hom_*(T, \Omega G), \beta \big)
\ar[dl]_{\iota_*} \\
& ~\pi_i\big(\Hom_*(T,G^I), \beta\big) \ar[r]^{\varphi_*} &
\pi_i\big(\Hom_*(T,G), \gamma\big) \ar[r] & \dots }$
\end{center}

From Lemma \ref{pathcontract}, we have that $\Hom_*(T,G^I)$ is
contractible, so that $\pi_i\big(\Hom_*(T,G^I)\big)=0$ for all $i$.
Hence the $\delta$ maps are all isomorphisms, and we get
$\pi_{i}\big(\Hom_*(T,G),\gamma\big) \cong \pi_{i-1}\big(\Hom_*(T,
\Omega G), \beta\big)$.  Applying this isomorphism $i-1$ times, we get

\begin{center}
$\pi_i\big(\Hom_*(T,G), \gamma\big) \cong \pi_1\big(\Hom_*(T,
\Omega^{i-1}(G)), \tilde\beta\big) \cong \pi_0 \big(\Hom_*(T,\Omega^i(G)\big) \cong [T, \Omega^i(G)]_\times$,

\end{center}

\noindent where the last two isomorphisms are from Lemma \ref{fundgroup} and Lemma
\ref{pointcomponent}, respectively.  This proves the claim.
\end{proof}

Since for a (pointed) topological space $\pi_i(X, x) =
\pi_{i-1}(\Omega X, \tilde x)$, and also $\pi_0\big(\Hom_*(T,G)\big)
= [T,G]_\times$, we see that in some sense the loop space functor
$\Omega$ commutes with the $\Hom_*$ complex, where it becomes the
graph theoretic version within the arguments.

\section{Concluding remarks}

With $\Omega G$ as our (pointed) loop space associated to a graph
$G$, we can define a graph theoretic notion of `$\times$-homotopy'
groups analogous to the $A$-homotopy groups from \cite{BBLL} (as
discussed in Section 2). For a pointed graph $G = (G,x)$, one can
interpret the path connected components of the graph $\Omega^n(G)$
as $\times$-homotopy classes of maps from the ``$n$ cube'' $I_m
\times I_m \dots \times I_m$ ($n$ times) into the graph $G$ such
that the boundary is mapped to the pointed vertex $x \in G$.  This
set is naturally a group under `stacking' and, by the above result,
is isomorphic to the group $\pi_n(\Hom_*({\bf 1}_*, G))$.  The
latter space is the clique complex on the subgraph of $G$ induced by
the looped vertices.

In the more general setting, we have a natural notion of the
`$T$-homotopy groups' of a pointed graph $G = (G,x)$ (for a fixed
finite pointed graph $T$).  These are defined according to
$\pi_i^T(G,x) := [T, \Omega^i(G)]_\times$; hence the groups described
in the previous paragraphs are obtained by setting $T = {\bf 1}_*$.
Note that this definition makes sense in any category `with a path
object' (see \cite{Bau} for an in depth discussion). Our results
show that $\pi_i^T(G,x) \cong \pi_i\big(\Hom_*(T,G)\big)$.

We have chosen to restrict our attention to ${\mathcal G}_*$, the category of pointed graphs.  This is the natural category to work in when one considers homotopy groups of spaces, which are by definition homotopy classes of certain pointed maps.  If one chooses to work in the category of (unpointed) graphs, the statement (and proof) of Theorem \ref{sequence} proceeds almost unchanged, and we obtain a long exact sequence of the relevant homotopy groups of (unpointed) $\Hom$ complexes.  
\begin{center}
$\xymatrix{ \dots \ar[r] & \pi_{i+1}\big(\Hom(T,G),\gamma\big)
\ar[r]^\delta & \pi_i\big(\varphi_T^{-1}(\leq \gamma), \beta\big)
\ar[dl]_{\iota_*}  \\
& ~\pi_i\big(\Hom(T,G^I), \beta\big)
\ar[r]^{\varphi_*} & \pi_i\big(\Hom(T,G), \gamma \big) \ar[r] &
\dots }$
\end{center}
One can check that in this case $\Hom(T,G^I) \simeq \Hom(T,G)$, as expected.  However, the $\big(\varphi_T^{-1}(\leq \gamma), \beta \big)$ term no longer has a natural interpretation as a $\Hom$ complex (as in the proof of Corollary \ref{Homgroups}).  In focussing on the pointed situation, we can apply the long exact sequence to show that the $\Hom_*$ complexes are spaces whose (usual) homotopy groups recover the combinatorically defined homotopy groups obtained from $\times$-homotopy of graph maps.

In conclusion, we wish to point out a certain similarity between our constructions and those of Schultz from \cite{Schspace}.  There it is shown that one can recover the space of (equivariant) maps from the circle into $\Hom(K_2, G)$ as a certain colimit (union) of $\Hom(C_{2r+1},G)$ complexes, where $C_{2r+1}$ is the odd cycle of length $2r + 1$.  In the process of proving this, the author implicitly uses the fact that the graph $C_{2r+1}^{K_2}$ is isomorphic to the interval graph $I_{4r+2}$ with the endpoints identified.  This suggests a general approach to understanding the topology of the mapping spaces of $\Hom$ complexes which is further explored in \cite{DS}.

\bibliographystyle{halpha.bst}
\bibliography{litgraph}

\begin{thebibliography}{BBdLL06}

\bibitem[Bau89]{Bau}
Hans~Joachim Baues.
\newblock {\em Algebraic homotopy}, volume~15 of {\em Cambridge Studies in
  Advanced Mathematics}.
\newblock Cambridge University Press, Cambridge, 1989.

\bibitem[BBdLL06]{BBLL}
Eric Babson, H{\'e}l{\`e}ne Barcelo, Mark de~Longueville, and Reinhard
  Laubenbacher.
\newblock Homotopy theory of graphs.
\newblock {\em J. Algebraic Combin.}, 24(1):31--44, 2006.

\bibitem[Bj{\"o}95]{Bjo95}
A.~Bj{\"o}rner.
\newblock Topological methods.
\newblock In {\em Handbook of combinatorics, Vol.\ 1,\ 2}, pages 1819--1872.
  Elsevier, Amsterdam, 1995.

\bibitem[BK06]{BKcom}
Eric Babson and Dmitry~N. Kozlov.
\newblock Complexes of graph homomorphisms.
\newblock {\em Israel J. Math.}, 152:285--312, 2006.

\bibitem[BK07]{BKpro}
Eric Babson and Dmitry~N. Kozlov.
\newblock {Proof of the Lov\'{a}sz Conjecture}.
\newblock {\em Annals of Mathematics}, 165(3):965--1007, 2007.

\bibitem[BL05]{BL}
H{\'e}l{\`e}ne Barcelo and Reinhard Laubenbacher.
\newblock Perspectives on {$A$}-homotopy theory and its applications.
\newblock {\em Discrete Math.}, 298(1-3):39--61, 2005.

\bibitem[Doca]{DocHom}
Anton Dochtermann.
\newblock {Hom complexes and homotopy theory in the category of graphs}.
\newblock to appear in \textit{European J. Combin.}

\bibitem[Docb]{DocUni}
Anton Dochtermann.
\newblock {The universality of Hom complexes}.
\newblock arXiv:math.CO/0702471.

\bibitem[DS]{DS}
Anton Dochtermann and Carsten Schultz.
\newblock {Constructing test graphs for topological bounds on chromatic
  number}.
\newblock Work in progress.

\bibitem[HN04]{HN04}
Pavol Hell and Jaroslav Ne{\v{s}}et{\v{r}}il.
\newblock {\em Graphs and homomorphisms}, volume~28 of {\em Oxford Lecture
  Series in Mathematics and its Applications}.
\newblock Oxford University Press, Oxford, 2004.

\bibitem[Kah07]{Ka}
Matthew Kahle.
\newblock {The neighborhood complex of a random graph}.
\newblock {\em J. Combin. Theory Ser. A}, 114(2):380--387, 2007.

\bibitem[Koz06]{Kfold}
Dmitry~N. Kozlov.
\newblock A simple proof for folds on both sides in complexes of graph
  homomorphisms.
\newblock {\em Proc. Amer. Math. Soc.}, 134(5):1265--1270 (electronic), 2006.

\bibitem[Koz07]{Kchr}
Dmitry~N. Kozlov.
\newblock Chromatic numbers, morphism complexes, and {S}tiefel-{W}hitney
  characteristic classes.
\newblock In {\em Geometric combinatorics}, volume~13 of {\em IAS/Park City
  Math. Ser.}, pages 249--315. Amer. Math. Soc., Providence, RI, 2007.

\bibitem[Lov78]{Lov78}
L.~Lov{\'a}sz.
\newblock Kneser's conjecture, chromatic number, and homotopy.
\newblock {\em J. Combin. Theory Ser. A}, 25(3):319--324, 1978.

\bibitem[ML98]{Mac98}
Saunders Mac~Lane.
\newblock {\em Categories for the working mathematician}, volume~5 of {\em
  Graduate Texts in Mathematics}.
\newblock Springer-Verlag, New York, second edition, 1998.

\bibitem[MZ04]{MZ}
Ji{\v{r}}{\'{\i}} Matou{\v{s}}ek and G{\"u}nter~M. Ziegler.
\newblock Topological lower bounds for the chromatic number: a hierarchy.
\newblock {\em Jahresber. Deutsch. Math.-Verein.}, 106(2):71--90, 2004.

\bibitem[Qui73]{Qui73}
Daniel Quillen.
\newblock Higher algebraic {$K$}-theory. {I}.
\newblock In {\em Algebraic $K$-theory, I: Higher $K$-theories (Proc. Conf.,
  Battelle Memorial Inst., Seattle, Wash., 1972)}, pages 85--147. Lecture Notes
  in Math., Vol. 341. Springer, Berlin, 1973.

\bibitem[Sch]{Schspace}
Carsten Schultz.
\newblock {Graph colourings, spaces of edges and spaces of circuits}.
\newblock To appear in \textit{Advances in Mathematics}.

\bibitem[{\v{Z}}iv05]{Ziv05}
Rade~T. {\v{Z}}ivaljevi{\'c}.
\newblock W{I}-posets, graph complexes and {${\Bbb Z}\sp 2$}-equivalences.
\newblock {\em J. Combin. Theory Ser. A}, 111(2):204--223, 2005.

\end{thebibliography}

\end{document}